\documentclass [12pt,leqno]{article}
\usepackage {latexsym}

\newcommand {\newproof}[2]%
  {\newenvironment{#1}%
  {\begin{trivlist}\item[\hskip \labelsep{\bf #2}]}%
  {\qed\end{trivlist}}}

\newcommand {\qed}{\hspace*{\fill}$\Box$}

\newcounter {claims}
\newenvironment {claims}
  {\begin{list}{\theclaims\hfill}%
  {\usecounter{claims}
  \renewcommand{\theclaims}{{\rm({\it\alph{claims}}\/)}}
  \labelwidth 2em
  \labelsep 0.5em
  \leftmargin 3.5em
  \rightmargin 0em}}%
  {\end{list}}

\newcounter {conditions}
\newenvironment {conditions}
  {\begin{list}{\theconditions\hfill}%
  {\usecounter{conditions}
  \renewcommand{\theconditions}{{\rm({\roman{conditions}})}}
  \labelwidth 2em
  \labelsep 0.5em
  \leftmargin 3.5em
  \rightmargin 0em}}%
  {\end{list}}

\newenvironment {properties}
  {\begin{conditions}}%
  {\end{conditions}}

\newcounter {savedconditions}
  {\setcounter {savedconditions}{\value{conditions}}
  \begin {conditions}
  \setcounter {conditions}{\value{savedconditions}}}%
  {\end{conditions}}

  {\setcounter {savedconditions}{\value{conditions}}
  \begin {properties}
  \setcounter {conditions}{\value{savedconditions}}}%
  {\end{properties}}

\newenvironment {conditionseq}
  {\begin{list}{\theconditions\hfill}%
  {\usecounter{conditions}
  \setcounter {conditions}{\value{equation}}
  \renewcommand{\theconditions}{{\rm({\arabic{conditions}})}}
  \labelwidth 2em
  \labelsep 0.5em
  \leftmargin 3.5em
  \rightmargin 0em}}%
  {\setcounter {equation}{\value{conditions}}
  \end{list}}

\newtheorem {theorem}{Theorem}
\newtheorem {lemma}{Lemma}
\newtheorem {corollary}{Corollary}
\newtheorem {proposition}{Proposition}

\newproof {proof}{\it Proof.}

\input {amssym}

\newcommand {\cf}{\mathop{\rm cf}\nolimits}
\newcommand {\ot}{\mathop{\rm ot}\nolimits}
\newcommand {\hgt}{\mathop{\rm ht}\nolimits}
\newcommand {\dom}{\mathop{\rm dom}\nolimits}
\newcommand {\ran}{\mathop{\rm ran}\nolimits}

\newcommand {\acc}{\mathop{\rm acc}\nolimits} 
\newcommand {\nacc}{\mathop{\rm nacc}\nolimits} 
\newcommand {\drop}{\mathop{\rm drop}\nolimits}
\newcommand {\gl}{\mathop{\mbox{$\rm g\ell$}}}

\newcommand {\GCH}{\mbox{\rm GCH}}

\newcommand {\ccc}{\mbox{\rm ccc}}
\newcommand {\NS}{\mbox{\rm NS}} 

\renewcommand {\emptyset}{\varnothing}
\newcommand {\sub}{\subseteq} 
\newcommand {\rest}{\mbox{$\mid$}} 

\newcommand {\pair}[2]{(#1,#2)}
\newcommand {\seq}[2]{(#1:#2)}
\newcommand {\singseq}[1]{(#1)}
\newcommand {\model}[1]{\langle#1\rangle}

\newcommand {\powerset}[1]{{\cal P}(#1)}

\newcommand {\card}[1]{{\left|#1\right|}}

\newcommand {\sing}[1]{\{#1\}} 
\newcommand {\set}[2]{\sing{#1:#2}}

\newcommand {\quine}[1]{\mbox{``#1''}}

\newcommand {\cand}{,\;} 

\newcommand {\cat}{\,^\frown} 

\newcommand {\union}{\cup}
\newcommand {\Union}{\bigcup}
\newcommand {\cut}{\cap} 
\newcommand {\Cut}{\bigcap}
\newcommand {\but}{\setminus} 

\newcommand {\leftscript}[2]{\,^{#1}#2}

\newcommand {\forces}{\Vdash}

\newcommand {\compat}{\Vert} 
\newcommand {\weak}[1]{{<#1}}

\newcommand {\pow}{\powerset}

\newcommand {\ls}{\leftscript}

\newcommand {\cc}[1]{\mbox{\rm $#1$-cc}} 

\newcommand {\mand}{\mbox{ and }}

\newcommand {\mforall}{\mbox{ for all }}

\newcommand {\eqtree}{7.13}

\title {Potential isomorphism and semi-proper trees}

\author {Alex Hellsten \thanks {Research partially supported by the
Mittag-Leffler Institute and grant 73499 from the Academy of Finland.
The authors wish to thank Pauli V\"ais\"anen for many helpful
suggestions.}
\and Tapani Hyttinen $^\ast$
\and Saharon Shelah \thanks {Research partially supported by the United
States-Israel Binational Science Foundation, the Israel Science
Foundation, and the Mittag-Leffler Institute.  Publication~770.}}
\date {September 11, 2001}

\begin {document}
\maketitle
\begin {abstract}
We study a notion of potential isomorphism, where two structures are
said to be potentially isomorphic if they are isomorphic in some
generic extension that preserves stationary sets and does not add new
sets of cardinality less than the cardinality of the models.  We
introduce the notions of semi-proper and weakly semi-proper trees, and
note that there is a strong connection between the existence of
potentially isomorphic models for a given complete theory and the
existence of weakly semi-proper trees.

We prove the existence of semi-proper trees under certain cardinal
arithmetic assumptions.  We also show the consistency of the
non-existence of weakly semi-proper trees assuming the consistency of
some large cardinals.
\end {abstract}

\section* {Introduction}

Two structures are said to be potentially isomorphic if they are
isomorphic in some extension of the universe in which they reside.
Different notions of potential isomorphism arise as restrictions are
placed on the method to extend the universe.  Nadel and
Stavi~\cite{NaSt} considered generic extensions in which there are no
new subsets of cardinality less than $\kappa$, where $\kappa$ is the
cardinality of the models.  They used some cardinal arithmetic
assumptions on $\kappa$ to show the existence of a pair of
non-isomorphic but potentially isomorphic models.  This kind of result
can be interpreted as a non-structure theorem for the theory of the
models in question.

In~\cite{HuHyRa:PI} these studies were continued, with an emphasis on
classification theory.  One of the results obtained there concerning
the notion introduced in~\cite{NaSt} is:
\begin {theorem} \label {unclassif} Let $T$ be a countable first order
theory and let \(\kappa = \kappa^{\aleph_0}\) be a regular cardinal.
The theory $T$ is unclassifiable if and only if there exists a pair of
non-isomorphic but potentially isomorphic models of $T$ of cardinality
$\kappa^+$.
\end {theorem}
A theory is said to be unclassifiable if it is unsuperstable or has
either the dimensional order property (DOP) or the omitting types
order property (OTOP).

Baldwin, Laskowski, and Shelah~\cite{BLSh:464,LwSh:518} studied a
weaker notion by considering extension by $\ccc$ forcing notions.
They showed that even classifiable theories may have a pair of
non-isomorphic models that are potentially isomorphic in this weaker
sense.

We must have some restrictions on how cardinals can be collapsed in
the extensions, because otherwise potential isomorphism will be
reduced to $L_{\infty \omega}$-equivalence.  But one may consider
weakening the requirement that the extension must be generic.  Such
notions are studied in~\cite{FrHyRa}, and it is shown there that this
kind of notions are not always decidable.  By a cardinal preserving
extension of $L$ we mean a transitive model of ZFC that contains all
ordinals, is contained in a set-generic extension of $V$, and has the
same cardinals as $L$.  For a tree \(T \in L\) on \({(\omega_1)}^L\),
let $C_T$ denote the set of all the trees \(T' \in L\) on
\({(\omega_1)}^L\) that are isomorphic with $T$ in some cardinal
preserving extension of $L$.  The following was proved
in~\cite{FrHyRa}:
\begin {theorem} Assume $0^\sharp$ exists.  There exists a tree \(T \in
L\) on \({(\omega_1)}^L\) such that $C_T$ is equiconstructible with
$0^\sharp$. \end {theorem}

The topic of this paper is a very strong notion of potential
isomorphism.  We consider generic extensions that preserve stationary
subsets of the cardinality of the models and do not add new sets of
cardinality less than the cardinality of the models.  To investigate
this notion of potential isomorphism is natural since
Theorem~\ref{unclassif} was proved in~\cite{HuHyRa:PI} by coding a
stationary set $S$ into a pair of models, that are then forced
isomorphic by killing $S$.

A \((\lambda, \kappa)\)-tree is a tree with the properties that every
branch has length less than $\kappa$ and every element has less than
$\lambda$ immediate successors.  Thus a \((\lambda, \kappa)\)-tree has
height at most $\kappa$.  Bearing some of the forthcoming proofs in
mind it is worth noting that the cardinality of a \((\lambda^+,
\kappa)\)-tree is at most $\lambda^{\weak{\kappa}}$.

We say that a \((\lambda, \kappa)\)-tree $T$ is \emph{weakly
semi-proper} if there exists a forcing notion $P$ that adds a
$\kappa$-branch to $T$, but preserves stationary subsets of $\kappa$
and adds no sets of cardinality less than $\kappa$.  If $T$ itself,
regarded as a forcing notion, has the properties of $P$ mentioned
above, then we say that $T$ is \emph{strongly semi-proper} or just
\emph{semi-proper}.

The following fact has lead us to questions concerning the existence
of weakly semi-proper \((\kappa^+, \kappa)\)-trees (for simplicity we
consider only countable theories):

\begin {theorem} Assume that $\kappa$ is uncountable and
\(\kappa^{\weak{\kappa}} = \kappa\).  The following statements are
equivalent:
\begin {conditions}
\item \label {SPT} There exists a weakly semi-proper \((\kappa^+,
\kappa)\)-tree
\item \label {SPTany} There exists a pair of non-isomorphic structures
of size $\kappa$ that can be made isomorphic by forcing, without
adding new sets of cardinality less than $\kappa$ or destroying
stationary subsets of $\kappa$
\item \label {SPTth} Statement \ref{SPTany} strengthened with the
requirement that the structures can be chosen to be models of any
complete countable theory $T$ such that either
\begin {enumerate}
\item $T$ is unstable,
\item $T$ has DOP, \(\kappa > {(c_r)}^+\), and
\(\xi^{c_r} < \kappa\) for every \(\xi < \kappa\), where $c_r$
is the smallest regular cardinal not less than the continuum, or
\item $T$ is superstable with DOP or OTOP.
\end {enumerate}
\end {conditions}
\end {theorem}

\begin {proof} \ref{SPTany} implies \ref{SPT}.  Suppose that two
non-isomorphic structures \(\goth A\) and \(\goth B\) of size $\kappa$
can be forced to be isomorphic without killing stationary sets or
adding new subsets of cardinality less than $\kappa$.  Let us assume
that $\kappa$ is the universe of both structures.  Let $P$ denote the
set of partial isomorphisms from \(\goth A\) to \(\goth B\) of
cardinality less than $\kappa$.  Let $T_\alpha$ denote the set
\[\set{f \in P}{\alpha \sub \dom f \cut \ran f \cand f [\kappa \but
\alpha] \union f^{-1} [\kappa \but \alpha] \sub \alpha}\] and let \(T
= \Union_{\alpha < \kappa} T_\alpha\) ordered by inclusion.  We shall
prove that $T$ is a \((\kappa^+, \kappa)\)-tree and that any forcing
notion that makes \(\goth A\) and \(\goth B\) isomorphic
without adding bounded subsets of $\kappa$ adds a $\kappa$-branch to
$T$.

It is straightforward to check that $T$ indeed is a tree.  Since
\(\kappa^\weak{\kappa} = \kappa\), the cardinality of $P$ is $\kappa$.
Therefore every node in $T$ has at most $\kappa$ immediate successors.
The union of a $\kappa$-branch would clearly be an isomorphism, so $T$
can not have $\kappa$-branches.  Finally suppose that $f$ is an
isomorphism between \(\goth A\) and \(\goth B\) in a generic
extension.  If there are no new bounded subsets of $\kappa$ in the
extension, then the function \((f \rest \alpha) \union (f^{-1} \rest
\alpha)^{-1}\) is in $T_\alpha$ for every \(\alpha < \kappa\) and it
follows that \(\pow{f} \cut T\) is a $\kappa$-branch through $T$ in
the generic extension.

\ref{SPT} implies \ref{SPTth}.  The proof of Lemma~\eqtree\ of
\cite{HuHyRa:PI} is essentially the proof of this implication.  It
relies on results of \cite{HyTu} and \cite{HySh:602}.
\end {proof}

Souslin trees are semi-proper \((\aleph_2, \aleph_1)\)-trees, and are
in fact used in that role in the proof of Lemma~\eqtree\ of
\cite{HuHyRa:PI}, but in this paper we shall see that semi-proper
trees exist under much weaker assumptions than Souslin trees.  The
following theorem summarises the rest of the results of this paper
except for some minor observations and strengthenings.
\begin {theorem} \label {MR} \hspace {1em}
\begin {claims}
\item \label {MRSC} It is consistent relative to a supercompact
cardinal that there are no weakly semi-proper \((\infty,
\aleph_1)\)-trees.
\item \label {MRCH} (Gregory) If \(2^{\aleph_0} < 2^{\aleph_1}\) then
there exists a semi-proper \((\aleph_2, \aleph_1)\)-tree.
\item \label {MRWC} It is consistent relative to a weakly compact
cardinal that there are no weakly semi-proper \((\aleph_3,
\aleph_2)\)-trees.
\item \label {MRGCH} Under $\GCH$ there exists a semi-proper
\((\kappa^{++}, \kappa^+)\)-tree for every infinite successor cardinal
$\kappa$.
\item \label {MR2} For any regular \(\kappa > \aleph_1\) there exists
a semi-proper \(({(2^\kappa)}^+, \kappa)\)-tree.
\end {claims}
\end {theorem}
Clause~\ref{MRWC} is proved in Section~\ref{collapsing}, \ref{MRCH}
and \ref{MRGCH} are proved in Section~\ref{diamonds}, and
Section~\ref{guessing} constitutes the proof of \ref{MR2}.
Clause~\ref{MRSC} follows from the observation that under Martin's
maximum (the semi-proper forcing axiom) there exists no weakly
semi-proper \((\infty, \aleph_1)\)-trees.  Feng~\cite{Feng:BP} has
made a similar observation concerning semi-proper \((\infty,
\aleph_1)\)-trees.

\section {Preliminaries and notation} \label {preli}

Let $A$ be a set of ordinals.  The set of ordinals $\alpha$ such that
\(\sup (A \cut \alpha) = \alpha\) (the accumulation points of $A$) is
denoted \(\acc^+ A\) and \(\acc A = \acc^+ A \cut A\) and \(\nacc A =
A \but \acc A\).  For infinite cardinals $\kappa$ and $\mu$ we let
$S^\kappa_\mu$ denote the set \(\set{\alpha \in \acc \kappa}{\cf
\alpha = \lambda}\).  $\NS_\kappa$ denotes the ideal of non-stationary
subsets of $\kappa$.  We sometimes use constructs like \(\alpha
(\ast)\) as ordinary variable names.  For sets $u$ and $E$ of ordinals
\[\drop (u,E) = \set{\sup (E \cut \alpha)}{\alpha \in u \cand \alpha >
\min E}.\] One can think of \(\drop (u,E)\) as the result of
``dropping'' $u$ onto $E$.  (In \cite{Sh:g} \(\drop (u,E)\) is denoted
\(\gl (u,E)\) where $\gl$ stands for ``glue''.)  Some of the
fundamental properties of $\drop$ that are needed in
Section~\ref{guessing} can be summarised as follows: If $E$ is closed
then \(\drop (u,E) \sub E\).  If $u$ is a club subset of some limit
ordinal $\delta$ and \(E \cut \delta\) is club in $\delta$ then
\(\drop (u,E)\) is club in $\delta$ and \(\acc (\drop (u,E)) \sub \acc
u \cut \acc E\).

In forcing arguments we follow the convention that \(p \leq q\) means
``$p$ is stronger than $q$''.  Our upward growing trees get inverted,
often without explicit mention, as soon as forcing with the tree in
question is discussed.

In Section~\ref{diamonds} we shall freely use some of the results
presented in~\cite{HuHyRa} about the ideal $I[\lambda]$ and the
$\kappa$-club game on a subset of $\lambda$, although we shall not
always stick to the notation used there.  The $\kappa$-club game on
\(S \sub \lambda\) is played by players I and II as follows: The game
lasts for $\kappa$ rounds.  On round $\xi$ player I first picks an
ordinal \(\alpha_\xi < \lambda\) that is greater than all the ordinals
played on earlier rounds.  Then player II picks an ordinal $\beta_\xi$
such that \(\alpha_\xi < \beta_\xi < \lambda\).  If the supremum of
the ordinals picked during the entire game is an element of $S$, then
player II wins the game.  Otherwise player I wins the game.  The game
characterisation of the $\kappa$-club filter on $\lambda$ is the
following statement:  If player II has a winning strategy in the
$\kappa$-club game on \(S \sub \lambda\) then there exists a set \(C
\sub S\) which is $\kappa$-club in $\lambda$.

\section {A consistency result} \label {collapsing}

We say that a tree $T$ is an \emph{$\alpha$-representation} (of a
tree) if the domain of $T$ is the ordinal $\alpha$ and \(x <_T y\)
implies \(x < y\) for all \(x,y \in T\).  Note that under the
assumption \(\kappa^\weak{\kappa} = \kappa\), every \((\kappa^+,
\kappa)\)-tree of height $\kappa$ is isomorphic to a
$\kappa$-representation.

\begin {lemma} \label {stat} If $\kappa$ is a regular uncountable
cardinal, $T$ is a $\kappa$-representation of a \((\kappa^+,
\kappa)\)-tree and the set \[S = \set{\alpha < \kappa}{\mbox{\(T \cut
\alpha\) has no $\alpha$-branch}}\] is stationary, then $T$ is not
weakly semi-proper.
\end {lemma}

\begin {proof} Suppose that $P$ is a forcing notion and $\dot{B}$ is a
$P$-name for a $\kappa$-branch through $T$.  Let $\dot{C}$ be a
$P$-name that satisfies
\[\forces \quine{\(\dot{C} = \set{\alpha < \kappa}{\ot (\dot{B}
\cut \alpha) = \alpha}\)}.\]
Assuming that $\kappa$ remains regular in the generic extension by
$P$, we get
\[\forces \quine{$\dot{C}$ is club and \(\dot{C} \cut \check{S} =
\emptyset\)}.\]
Thus $P$ necessarily kills a stationary set, which shows that $T$ can
not be weakly semi-proper.
\end {proof}

\begin {corollary} \label {WC}
If $\kappa$ is weakly compact then weakly semi-proper \((\kappa^+,
\kappa)\)-trees do not exist.
\end {corollary}

\begin {proof} Let $T$ be a $\kappa$-representation for a \((\kappa^+,
\kappa)\)-tree.  The fact that $T$ has no cofinal branches can be
expressed as a $\Pi^1_1$-statement in the structure \(\model{V_\kappa,
\in, T}\).  For regular \(\alpha < \kappa\) the same
$\Pi^1_1$-statement interpreted in \(\model{V_\alpha, \in, T \cut
V_\alpha}\) expresses the fact that \(T \cut \alpha\) has no
$\alpha$-branches.  Given this $\Pi^1_1$-statement, the corollary
immediately follows from Lemma~\ref{stat} by
$\Pi^1_1$-reflection. \end {proof}

We shall now give the definition of a forcing notion that was
introduced by Mitchell~\cite{Mi}.  Let $\kappa$ be a weakly compact
cardinal.  Let $P$ be the classical forcing notion for adding $\kappa$
many Cohen reals.  In other words $P$ is the set of finite partial
functions from $\kappa$ to $2$, ordered by reverse inclusion.  Let
$B(P)$ be the complete boolean algebra associated with $P$.  For \(s
\sub P\) we shall use the notation $b_s$ for the regular open cover
(see e.g.\ {Jech}~\cite[Lemma~17.2]{Jech}) of $s$, so that we have
\(B(P) = \set{b_s}{s \sub P}\).

Let \(P_\alpha = \set{p \in P}{p \rest \alpha = p}\) and \(B_\alpha =
\set{b_s}{s \sub P_\alpha}\).  Then $B_\alpha$ is isomorphic to
$B(P_\alpha)$.  A partial function \(f: \kappa \to B(P)\) is
acceptable if \(\card{f} < \aleph_1\) and \(f(\gamma) \in B_{\gamma +
\omega}\) for every \(\gamma < \kappa\).  We let $A$ denote the set of
all acceptable functions.  Given a $P$-generic set $G$, we define a
forcing notion $Q$ in $V[G]$ as follows: For every \(f \in A\), where
$A$ is regarded as an element of $V$, let $\bar{f}$ denote the
characteristic function of \(\set{\gamma \in \dom f}{f(\gamma) \cut G
\neq \emptyset}\).  Then let $Q$ be \(\set{\bar{f}}{f \in A}\) ordered
by reverse inclusion.  With $\dot{Q}$ being a $P$-name for $Q$, we
finally let $R$ be the two step iteration \(P \ast \dot{Q}\).  We
shall also refer to $R$ as the Mitchell forcing.  The model $V^R$
obtained by assuming $\GCH$ and then forcing with $R$, we shall call
the Mitchell model.  In the notation of~\cite{Mi} our $R$ is
isomorphic to \(R_2 (\aleph_0, \aleph_1, \kappa)\).

Let \(Q_\alpha = \set{\bar{f} \in Q}{\bar{f} \rest \alpha = \bar{f}}\)
and let \(R_\alpha = P_\alpha \ast \dot{Q}_\alpha\) where the ordering
of $Q_\alpha$ is reverse inclusion.  Thus \(R_\kappa = R\).  For any
$R$-generic set $G$, we let $G_\alpha$ denote the set \(G \cut
R_\alpha\).  We shall need the following results from~\cite{Mi}:

\begin {lemma}[Mitchell] \label {prop} Assume that $\GCH$ holds.
\begin {claims}
\item Suppose that $\alpha$ is a limit ordinal in $\kappa$ and $G$ is
a $R$-generic set.  Then $G_\alpha$ is
$R_\alpha$-generic.
\item \label {all-in} Suppose that \(\cf \gamma > \omega\) and $f$ is
a function \(\gamma \to V\) in $V^R$.  If \(f \rest \zeta \in
V^{R_\alpha}\) for every \(\zeta < \gamma\) then \(f \in
V^{R_\alpha}\).
\item $R$ has the $\cc{\kappa}$.
\item \label {card} In $V^R$, \(2^{\aleph_1} = \kappa =
\aleph_2\).
\end {claims}
\end {lemma}

\begin {proposition} \label {Mitch}
In the Mitchell model there are no weakly semi-proper \((\aleph_3,
\aleph_2)\)-trees. \end {proposition}

\begin {proof} Let $R_\kappa$ be the Mitchell forcing notion and let
$\dot{T}$ be a $R_\kappa$-name for an arbitrary \((\aleph_3,
\aleph_2)\)-tree.  By clause~\ref{card} of Lemma~\ref{prop} we can
assume that $\dot{T}$ is a name for an $\omega_2$-representation and
by Lemma~\ref{stat} it is then enough to prove that
\[\forces_{R_\kappa} \quine{\(\set{\alpha <
\omega_2}{\mbox{\(\dot{T} \cut \alpha\) has no $\alpha$-branch}}\) is
stationary}.\] Since $R_\kappa$ is $\cc{\kappa}$ and therefore does
not destroy stationary sets, it is even sufficient to find a
stationary set \(S \sub \kappa\), such that
\begin {equation} \label {S}
\forces_{R_\kappa} \quine{\(\dot{T} \cut
\alpha\) has no $\alpha$-branch when \(\alpha \in \check{S}\)}.
\end {equation}

We shall use $\Pi^1_1$-reflection to find a stationary set $S$
satisfying (\ref{S}).  To be able to capture various facts about
forcing using $\Pi^1_1$-statements in a structure like
\(\model{V_\kappa, \in, R_\kappa, \dot{T}}\) we need to make some
assumptions about the names used.  The name $\dot{T}$ can be assumed
to be a subset of \((\kappa \times \kappa) \times R_\kappa\) where we
identify ordinals with their canonical names.  Furthermore we can
assume that for every \(\pair{\alpha}{\beta} \in \kappa \times
\kappa\) the set
\[A_{\pair{\alpha}{\beta}} = \set{p \in
R_\kappa}{\pair{\pair{\alpha}{\beta}}{p} \in \dot{T}}\] is a maximal
antichain of the set consisting of all conditions $p$ with the
property \(p \forces \pair{\alpha}{\beta} \in \dot{T}\).  Then for any
\(q \in R_\kappa\), \(q \forces \pair{\alpha}{\beta} \notin \dot{T}\)
if and only if \(\set{p \in A_{\pair{\alpha}{\beta}}}{p \compat q}\)
is empty.  An arbitrary name for a subset of $\dot{T}$ can be thought
of as a name for a subset of $\kappa$ and then there always exists an
equivalent name that is a subset of \(\kappa \times R_\kappa\) and has
similar properties as $\dot{T}$ above.  For such a name $\dot{B}$ for
a subset of $\dot{T}$ the statement \[\forces_{R_\kappa}
\quine{$\dot{B}$ is a $\kappa$-branch through $\dot{T}$}\] can be
expressed with a first order sentence in the structure
\(\model{V_\kappa, \in, R_\kappa, \dot{T}, \dot{B}}\).  Let us call a
name like $\dot{T}$ or $\dot{B}$ \emph{normal} for the rest of the
proof.  Normality of a name is also a first order property of the
structure mentioned above.

For inaccessible cardinals \(\alpha < \kappa\), \(R_\kappa \cut
V_\alpha = R_\alpha\) and if we let \(\dot{T}_\alpha = \dot{T} \cut
V_\alpha\) and \(\dot{B}_\alpha = \dot{B} \cut V_\alpha\) then
$\dot{T}_\alpha$ and $\dot{B}_\alpha$ are $R_\alpha$-names.  So there
is a $\Pi^1_1$-sentence $\sigma$ such that for every inaccessible
\(\alpha \leq \kappa\), \(\model{V_\alpha, \in, R_\alpha, T_\alpha}
\models \sigma\) if and only if $\dot{T}_\alpha$ is normal and
\begin {equation} \label {branch}
\forces_{R_\alpha} \quine{$\dot{T}_\alpha$ has no $\alpha$-branch.}
\end {equation}
Furthermore there exists a club subset $D$ of $\kappa$ such that
\begin {equation} \label {G}
(\dot{T}_\alpha)_{G_\alpha} = \dot{T}_G \cut \alpha
\end {equation}
for every \(\alpha \in D\) and every $R_\kappa$-generic set $G$.  Let
$S$ be a stationary set of ordinals such that (\ref{branch}) and
(\ref{G}) hold for every \(\alpha \in S\).  By clause~\ref{all-in} of
Lemma~\ref{prop} it now follows that $S$ satisfies (\ref{S}).
\end {proof}

\section {Using weak diamond principles} \label {diamonds}

We say that a tree $T$ is \emph{splitting} if it has unique limits and
if every node of $T$ has at least two immediate successors.  If $T$ is
splitting and for every \(x \in T\) and \(\alpha < \hgt T\) there
exist an element \(y \in T\) such that \(x <_T y\) and \(\hgt y \geq
\alpha\), then we say that $T$ is \emph{normal}.  Let $\kappa$ be
regular and uncountable and let $T$ be a normal tree of height
$\kappa$.  If forcing with $T$ adds a new set of cardinality less than
$\kappa$, then $\kappa$ becomes singular in the generic extension.
Thus if forcing with $T$ preserves stationary subsets of $\kappa$,
then no new sets of cardinality less than $\kappa$ are added.

A subset $U$ of a tree $T$ is called a {\em $\mu$-fan} of $T$ if there
exists a sequence \(\seq{\delta_\xi}{\xi < \mu}\) and an indexed
family \(\seq{x_f}{f \in \ls{\weak{\mu}}{2}}\) such that
\begin {conditionseq}
\item \label {fanA} \(U = \set{x_f}{f \in \ls{\weak{\mu}}{2}}\)
\item \(\seq{\delta_\xi}{\xi < \mu}\) is strictly increasing and continuous
\item \(\hgt_T x_f = \delta_{\dom f}\) for every \(f \in \ls{\weak{\mu}}{2}\)
\item \label {fanO} \(\inf_T \sing{x_{f \cat \singseq{0}}, x_{f \cat
\singseq{1}}} = f_x\) for every \(f \in \ls{\weak{\mu}}{2}\).
\end {conditionseq}
We say that $T$ is {\em $\mu$-fan closed} if $T$ is $\mu$-closed as a
forcing notion, and for every $\mu$-fan $U$ of $T$ there exists an
element \(x \in T\) that extends one of the cofinal branches in $U$.

\begin {lemma} \label {fan} Suppose that \(\mu^{\weak{\mu}} = \mu\)
and \(\kappa = \mu^+\).  Then every splitting $\mu$-fan closed
\((\infty, \kappa)\)-tree is semi-proper. \end {lemma}

\begin {proof} It is straight forward to prove by induction that a
splitting $\mu$-fan-closed \((\infty, \kappa)\)-tree must be a normal
tree of height $\kappa$.  By normality forcing with the tree must
produce a $\kappa$-branch.  Thus it only remains to prove that
stationary sets are preserved.

Let $P$ be an inverted normal $\mu$-fan closed tree of height
$\kappa$, let \(S \sub \kappa\), and let $\dot{C}$ be a $P$-name such
that \[\forces \quine{$\dot{C}$ is club and \(\dot{C} \cut \check{S} =
\emptyset\)}.\] Because $I[\kappa]$ is improper by our assumptions,
the game characterisation of the $\mu$-club filter on $\kappa$ holds.
We shall finish the proof by showing that player~II has a winning
strategy in the $\mu$-club game on the complement of $S$.  This will
be enough since we can assume that \(S \sub S^\kappa_\mu\).  The
strategy can be described as follows.  At round $\xi$ in the game,
player~I has picked $\alpha_\xi$ and player~II should now answer with
\(\beta_\xi > \alpha_\xi\).  But before fixing $\beta_\xi$ we pick a
set \(\set{p_f}{f \in \ls{\xi}{2}}\) of conditions in $P$ and a set
\(\set{\gamma_f}{f \in \ls{\xi}{2}}\) of ordinals such that the
following holds for every $f$ and $g$ in $\ls{\xi}{2}$:
\begin {conditionseq}
\item \(p_f \leq p_{f \rest \nu}\) for every \(\nu < \xi\)
\item If \(f \neq g\) then \(\hgt p_f = \hgt p_g\) and if \(\xi =
\nu+1\) then \(\sup \sing{p_f, p_g} = p_{f \rest \nu}\)
\item If $\xi$ is a limit ordinal then \(\hgt p_f = \sup_{\nu < \xi}
\hgt p_{f \rest \nu}\)
\item \(\hgt p_f > \gamma_h\) for every \(h \in \Union_{\nu < \xi}
\ls{\nu}{2}\)
\item \(\gamma_f > \alpha_\xi\) and if $\xi$ is a successor ordinal
then \(p_f \forces \gamma_f \in \dot{C}\).
\end {conditionseq}
Then we put \(\beta_\xi = \sup \set{\gamma_f \union \hgt p_f}{f \in
\ls{\xi}{2}}\) if $\xi$ is a successor ordinal and \(\beta_\xi =
\alpha_\xi + 1\) otherwise.  Let \(\alpha = \sup_{\xi < \mu}
\alpha_\xi\).  Since \(\set{p_f}{f \in \Union_{\xi < \mu}
\ls{\xi}{2}}\) is a $\mu$-fan, there exists a function \(f: \mu \to
2\) and a condition $p$ such that \(p \leq p_{f \rest \xi}\) for every
\(\xi < \mu\).  Now \(p \forces \alpha \in \dot{C}\) which implies
that \(\alpha \notin S\).
\end {proof}

The combinatorial principle called weak diamond defined in
\cite{DvSh:65} is equivalent with \(2^{\aleph_0} < 2^{\aleph_1}\).
The tree construction in the proof below is essentially due to
Gregory~\cite{Gr:BA}.  The proof is shortened considerably by the use
of the weak diamond principle of~\cite{DvSh:65} which is implicitly
proved in Gregory's construction.

\begin {proposition}[Gregory] If \(2^{\aleph_0} < 2^{\aleph_1}\)
then there exists a semi-proper \((\aleph_2, \aleph_1)\)-tree.
\end {proposition}

\begin {proof} We can recursively define a function \(F:
\ls{\weak{\omega_1}}{2} \to 2\) with the following property: Every
$\aleph_0$-fan of $\ls{\weak{\omega_1}}{2}$ has two cofinal branches
such that if $x$ and $y$ are the unions of these branches then \(F(x)
\neq F(y)\).  By the weak diamond principle there exists a function
\(g: \omega_1 \to 2\) such that \(\set{\alpha < \omega_1}{F(f \rest
\alpha) = g(\alpha)}\) is stationary for every \(f: \omega_1 \to 2\).
Clearly \[T = \set{f \in \ls{\weak{\omega_1}}{2}}{F(f \rest \alpha)
\neq g(\alpha) \mforall \alpha \in \acc^+ (\dom f)}\] is a splitting
\((\aleph_2, \aleph_1)\)-tree.  The function $F$ was constructed in
such a way that $T$ is guaranteed to be $\aleph_0$-fan closed.  Then
$T$ is a semi-proper \((\aleph_2, \aleph_1)\)-tree by
Lemma~\ref{fan}. \end {proof}

Let $E$ be a stationary subset of $\kappa^+$ where $\kappa$ is some
infinite cardinal.  For \(\delta \in E\), let \(\eta_\delta: \cf
\delta \to \delta\) be an increasing continuous function with limit
$\delta$.  We let \(\Phi (\eta_\delta: \delta \in E)\) denote the
following combinatorial principle: There exists a sequence
\(\seq{d_\delta}{\delta \in E}\) where each $d_\delta$ is a function
\(\cf \delta \to \delta\) such that for any function \(h: \kappa^+ \to
2\), there is a stationary set of ordinals \(\delta \in E\) satisfying
\[\set{i < \cf \delta}{d_\delta (i) = h (\eta_\delta (i))} \mbox{ is
stationary in } \cf \delta.\] The sequence \(\seq{d_\delta}{\delta \in
E}\) can be referred to as a {\em weak diamond} sequence.

We shall use the following result by Shelah~\cite[Appendix,
Theorem~3.6]{Sh:f}:
\begin {lemma} If \(\kappa = \kappa^{\weak{\kappa}}\) and \(\kappa =
2^\theta\) for some cardinal $\theta$, then \(\Phi (\eta_\delta:
\delta \in S^{\kappa^+}_\kappa)\) holds for any sequence
\(\seq{\eta_\delta}{\delta \in S^{\kappa^+}_\kappa}\) as defined
above. \end {lemma}

\begin {proposition} \label {plusplus} If \(\kappa = \theta^+ =
2^\theta\) for some cardinal $\theta$ then there exists a semi-proper
\((\kappa^{++}, \kappa^+)\)-tree. \end {proposition}

\begin {proof} Let \(E = S^{\kappa^+}_\kappa\), fix
\(\seq{\eta_\delta}{\delta \in E}\), and let \(\seq{d_\delta}{\delta
\in E}\) be a weak diamond sequence given by \(\Phi (\eta_\delta:
\delta \in E)\).  We claim that \[T = \set{f \in
\ls{\weak{\kappa^+}}{2}}{\forall \delta \in E \cut \acc^+ (\dom f)
(\set{i < \kappa}{d_\delta (i) = f (\eta_\delta (i))} \in
\NS_{\kappa})}\] is the required tree.  Clearly $T$ is a splitting
\((\kappa^{++}, \kappa^+)\)-tree.  By Lemma~\ref{fan} it then suffices
to prove that $T$ is $\kappa$-fan closed.

It is immediate from the definition that $T$ is $\kappa$-closed.  Let
$U$ a $\kappa$-fan of $T$ and suppose that \(\seq{x_f}{f \in
\ls{\weak{\kappa}}{2}}\) and the sequence \(\seq{\delta_\xi}{\xi <
\kappa}\) satisfy conditions \ref{fanA}--\ref{fanO}.  Let \(\delta =
\sup_{\xi < \kappa} \delta_\xi\).  By~\ref{fanO} we may assume without
loss of generality that
\[x_f (\delta_\xi) = f(\xi) \mforall \xi < \kappa \mand f: \xi+1 \to 2.\]
Now we make use of the fact that \(\set{\delta_\xi}{\xi < \kappa} \cut
\ran \eta_\delta\) is a club subset of $\delta$.  Define a function
\(f: \kappa \to 2\) by letting \(f (\nu) = 1 - d_\delta (i)\) whenever
\(\eta_\delta (i) = \delta_\nu\).  Now \(\Union_{\xi < \kappa} x_{f
\rest \xi}\) is in $T$ which shows that $T$ is $\kappa$-fan
closed. \end {proof}

\section {Semi-proper trees in ZFC} \label {guessing}

This entire section constitutes the proof of clause~\ref{MR2} of
Theorem~\ref{MR}.  For convenience we restate the result.

\begin {proposition} \label {part2}
For any regular \(\kappa > \aleph_1\) there exists a
semi-proper \(((2^\kappa)^+, \kappa)\)-tree. \end {proposition}
We shall first define a tree $T$ as a subtree of \(\Union_{\alpha <
\kappa} \ls{\alpha + 1}{\pow{\kappa}}\) ordered by inclusion.  $T$
will be a semi-proper \(({(2^\kappa)}^+, \kappa)\)-tree unless it has
a $\kappa$-branch.  If $T$ has a $\kappa$-branch we shall use this
branch to construct another tree that meets the requirements.  In fact
this second tree will be a semi-proper \((\kappa^+, \kappa)\)-tree.

For functions \(p: \alpha + 1 \to \pow{\kappa}\) we shall use the
following notation.  The ordinal $\alpha$ is denoted \(\alpha (p)\).
For every \(\beta \leq \alpha\)
\[u_\beta = \left \{ \begin{array}{ll}
p(\beta) & \mbox{if $p(\beta)$ is a closed subset of $\beta$} \\
\emptyset & \mbox{otherwise}
\end{array} \right .\]
and \[S_\beta = \left \{ \begin{array}{ll} p(\beta) & \mbox{if
$p(\beta)$ is stationary in $\kappa$} \\ \emptyset & \mbox{otherwise.}
\end{array} \right. \] We write $u^p_\beta$ and $S^p_\beta$ for
$u_\beta$ and $S_\beta$ respectively, if $p$ is not clear from the
context.

\subsection* {The first tree}

We let \(p \in T\) if and only if the following conditions
hold whenever \(\gamma < \beta \leq \alpha\):
\begin {conditionseq}
\item If $u_\beta$ is empty then $S_\beta$ is non-empty (thus
stationary)
\item \label {P2} If \(\gamma \in u_\beta\) then \(u_\gamma = u_\beta
\cut \gamma\)
\item \label {P3} If $\beta$ is a limit ordinal then $u_\beta$ is
unbounded in $\beta$
\item \label {P9} If \(\gamma \in u_\beta\) and $\gamma$ is a limit
then \(\gamma \notin S_{\min u_\gamma}\).
\end {conditionseq}

We shall now prove that forcing with $T$ does not destroy stationary
subsets of $\kappa$.  Let $S$ be a stationary set, let \(p \in T\),
and let $\dot{C}$ be a name that is forced by $p$ to be club in
$\kappa$.  We construct a condition \(q \leq p\) such that \(q \forces
\dot{C} \cut \check{S} \neq \emptyset\).  By induction on \(i <
\kappa\) we continue for as long as possible to pick conditions $p_i$
and ordinals $\alpha_i$ such that the following holds when $p_i$ and
$\alpha_i$ have been defined for every \(i < \zeta\):
\begin {conditionseq}
\item \label {US} \(S^{p_0}_{\alpha_0} = S\)
\item \label {UD}\(\seq{p_i}{i < \zeta}\) is decreasing and \(p_0 \leq
p\)
\item \label {UI}\(\seq{\alpha_i}{i < \zeta}\) is increasing and
continuous
\item \label {UC}\(p_{i+1} \forces \dot{C} \cut (\alpha_{i+1} \but
\alpha_i) \neq \emptyset\)
\item \label {UU}\(\alpha (p_i) \geq \alpha_i\) (alternatively
\(\alpha (p_i) = \alpha_i\)) and \(u^{p_i}_{\alpha_i} =
\set{\alpha_j}{j < i}\)
\item \label {UL}If $\alpha_i$ is a limit then $i$ is a limit and
\(\alpha_i \notin S\).
\end {conditionseq}
We shall drop the superscripts on $u^{p_i}_\beta$ and $S^{p_i}_\beta$
because condition~\ref{UD} makes them obsolete.  Clearly we can put
\(p_0 = p \cat \singseq{S}\) and \(\alpha_0 = \alpha (p) + 1\).  We
shall now check that appropriate $p_{i+1}$ and $\alpha_{i+1}$ always
can be picked once the preceding conditions and ordinals have been
successfully defined.  First pick \(q \leq p_i\) and \(\gamma \geq
\alpha_i\) such that \(q \forces \gamma \in \dot{C}\).  Then let
\(\alpha_{i+1} = \max \sing{\alpha (q), \gamma} + 1\).  Now we shall
define \(p_{i+1}: \alpha_{i+1} + 1 \to \pow{\kappa}\) by fixing
$u_\beta$ and $S_\beta$ for ordinals $\beta$ such that \(\alpha (q) <
\beta \leq \alpha_{i+1}\).  Let \(u_{\alpha_{i+1}} = \set{\alpha_j}{j
< i+1}\) and if \(\alpha_{i+1} > \alpha (q) + 1\), let \(S_{\alpha (q)
+ 1} = \kappa \but \alpha_{i+1}\).  Finally fill the possible gap by
letting \(u_\beta = \beta \but (\alpha (q) + 1)\) for those ordinals
$\beta$ that satisfy \(\alpha (q) + 1 < \beta < \alpha_{i+1}\).

Now suppose that we are about to pick $p_i$ where $i$ is a limit.
By~\ref{UI} we must have \(\alpha_i = \Union_{j < i} \alpha_j\) in
this situation.  The only possible way to define $p_i(\alpha_i)$ is to
let \(u_{\alpha_i} = \set{\alpha_j}{j < i}\).  Let \(q = (\Union_{j <
i} p_j) \cat \singseq{u_{\alpha_i}}\).  If $\alpha_i$ happens to be in
the complement of $S$, we can make the induction go on by putting
\(p_i = q\).  But if \(\alpha_i \in S\) we are done with the proof
because, in any case, \(q \forces \alpha_i \in \dot{C}\).  The latter
must happen sooner or later because otherwise we finally have \(S \cut
\acc \set{\alpha_i}{i < \kappa} = \emptyset\) contradicting the
assumption that $S$ is stationary.

The proof that $T$ is normal is similar to the successor step in the
construction above.  If $T$ does not have cofinal branches then the
proposition is proved.  Let us now assume that $T$ has a cofinal
branch and construct another tree that has the required properties.

\subsection* {The second tree}

The cofinal branch through $T$ gives us two sequences
\(\seq{u_\beta}{\beta < \kappa}\) and \(\seq{S_\beta}{\beta <
\kappa}\) such that $u_\beta$ is a closed subset of $\beta$ and
$S_\beta$ is stationary in $\kappa$ for every \(\beta < \kappa\) and
the conditions \ref{P2}--\ref{P9} hold.  For every \(\alpha < \kappa\)
let
\begin {equation} \label {Sstar}
S^\ast_\alpha = \set{\beta < \kappa}{\alpha \in u_\beta}
\end {equation}
and let $E_\alpha$ be a club subset of $\kappa$ such that
\(S^\ast_\alpha \cut E_\alpha = \emptyset\) whenever $S^\ast_\alpha$
is non-stationary.  Let $E$ be the diagonal intersection \(\set{\beta
< \kappa}{\beta \in \Cut_{\alpha < \beta} E_\alpha}\).  It is now easy
to verify that if \(\beta \in E\) then $S^\ast_\alpha$ is stationary
for every \(\alpha \in u_\beta\).

\begin {lemma} \label {ast} There exist ordinals \(\alpha (\ast)\)
and \(\beta (\ast)\) such that \(\alpha (\ast) < \beta (\ast) <
\kappa\), \(S^\ast_{\alpha (\ast)}\) and \(S^\ast_{\beta (\ast)} \cut
S^\kappa_{\omega_1}\) are stationary in $\kappa$, and \(S^\ast_{\alpha
(\ast)} \cut S^\ast_{\beta (\ast)} = \emptyset\).
\end {lemma}

\begin {proof}
First we shall find limit ordinals \(\alpha, \beta \in E\) such that
\(\alpha < \beta\) and \(\alpha \notin u_\beta\).  Let $\alpha$ be a
limit ordinal in $E$ and let \(\beta > \alpha\) be a limit ordinal in
\(E \cut S_{\min u_\alpha}\).  Let \(\gamma > \beta\) be limit ordinal
in $E$.  If \(\alpha \in u_\beta\) then \(\beta \notin u_\gamma\) so
the required ordinals can be picked by replacing, if necessary,
$\alpha$ and $\beta$ by $\beta$ and $\gamma$ respectively.

Fix \(\alpha (\ast) \in u_\alpha\) such that \(\alpha (\ast) > \Union
(u_\beta \cut \alpha)\) and let \(\beta (\ast) = \min (u_\beta \but
\alpha)\).  From what was noted above about $E$ it is now clear that
\(S^\ast_{\alpha (\ast)}\) and \(S^\ast_{\beta (\ast)}\) are both
stationary and disjoint from each other.  We shall now prove that
\(S^\ast_{\beta (\ast)} \cut S^\kappa_{\omega_1}\) can be assumed to
be stationary.  Suppose that $C$ is a club such that \(S^\ast_{\beta
(\ast)} \cut S^\kappa_{\omega_1} \cut C = \emptyset\).  Define a
function \(f: S^\kappa_{\omega_1} \but \beta (\ast) \to \kappa\) by
\(f (\gamma) = \min (u_\gamma \but \beta (\ast))\).  By Fodor's lemma
there exists a stationary set \(S \sub S^\kappa_{\omega_1} \cut C\)
and an ordinal \(\delta (\ast)\) such that \(f[S] = \sing{\delta
(\ast)}\).  Now \(S^\ast_{\delta (\ast)} \cut S^\kappa_{\omega_1}\) is
stationary because it has $S$ as a subset.  We must have \(\beta
(\ast) \notin u_{\delta (\ast)}\) and \(\beta (\ast) < \delta (\ast)\)
because \(\beta (\ast) \in u_{\delta (\ast)}\) or \(\beta (\ast) =
\delta (\ast)\) would imply that \(S^\ast_{\delta (\ast)} \sub
S^\ast_{\beta (\ast)}\) which contradicts the assumption that
\(S^\ast_{\beta (\ast)} \cut S^\kappa_{\omega_1}\) is non-stationary.
But this means that \(S^\ast_{\beta (\ast)}\) and \(S^\ast_{\delta
(\ast)}\) are disjoint and could thus serve as replacements for
\(S^\ast_{\alpha (\ast)}\) and \(S^\ast_{\beta (\ast)}\) respectively.
\end {proof}

Fix ordinals \(\alpha (\ast)\) and \(\beta (\ast)\) with the
properties stated in the last lemma.  Next we shall construct a ``club
guessing'' sequence that can be used in tree construction in a similar
way as the weak diamond principles presented in Section~\ref{diamonds}.

\begin {lemma} \label {CG} There exists a club \(E^\ast \sub \acc \kappa\)
and a sequence \(\seq{C_\delta}{\delta \in S^\ast_{\beta (\ast)} \cut
\acc E^\ast}\) such that
\begin {conditionseq}
\item \label {CGclub} $C_\delta$ is club in $\delta$
\item \label {CGnacc} \(C_\delta \cut S^\ast_{\alpha (\ast)} \sub
\nacc C_\delta\)
\item \label {CGstat} For any club \(E' \sub E^\ast\) the set
\[\set{\delta \in S^\ast_{\beta (\ast)} \cut \acc E^\ast}{\delta =
\sup (E' \cut \nacc C_\delta \cut S^\ast_{\alpha (\ast)})}\] is
stationary in $\kappa$
\item \label {CGcoh} \(\delta' \in u_{\delta} \cut S^\ast_{\beta
(\ast)} \cut \acc E^\ast\) implies \(C_{\delta'} = C_{\delta} \cut
\delta'\).
\end {conditionseq}
\end {lemma}

\begin {proof} 
Let \(E_0 = \acc \kappa\) and let \(C^0_\delta = \drop (u_\delta,
E_0)\) for every \(\delta \in S^\ast_{\beta (\ast)} \cut \acc E_0\).
By recursion on $n$ we define club sets $E_n$ and sequences
\(\seq{C^n_\delta}{\delta \in S^\ast_{\beta (\ast)} \cut \acc E_n}\)
such that \(E_{n+1} \sub \acc E_n\),
\begin {equation} \label {nextE}
\delta > \sup (E_{n+1} \cut \nacc C^n_\delta \cut S^\ast_{\alpha
(\ast)}) \mforall \delta \in S^\ast_{\beta (\ast)} \cut E_{n+1},
\end {equation}
and $C^{n+1}_\delta$ is defined by
\begin {equation} \label {nextC}
C^{n+1}_\delta = C^n_\delta \union \Union_\beta \drop
(u_\beta, E_{n+1}) \but \gamma^n_\delta (\beta)
\end {equation}
where the large union is taken over all \(\beta \in (\nacc
C^n_\delta) \but (S^\ast_{\alpha (\ast)} \cut E_{n+1})\) and
\begin {equation} \label {gamma}
\gamma^n_\delta (\beta) = \left\{ \begin{array}{ll} \max
((C^n_\delta \cut \beta) \union \sing{0}) & \mbox{if \(\sup (E_{n+1}
\cut \beta) = \beta\)} \\ \max ((E_{n+1} \cut \beta) \union \sing{0})
& \mbox{otherwise.} \end{array} \right.
\end {equation}
We claim that for some \(n < \omega\) there exists no club \(E_{n+1}
\sub E_n\) satisfying (\ref{nextE}), and that when this happens the
sets \(C_\delta = C^n_\delta\) and the set \(E^\ast = E_n\) satisfy
the conditions of the lemma.

In fact it is straightforward (see Section~\ref{preli}) to check that
conditions \ref{CGclub}, \ref{CGnacc}, and \ref{CGcoh} hold for every
\(n < \omega\) even if we drop the requirement (\ref{nextE}) and just
pick any club \(E_{n+1} \sub \acc E_n\) during the construction.  To
see by induction that \ref{CGclub} and \ref{CGnacc} hold, let
\(\seq{\alpha_i}{i < \zeta}\) be a strictly increasing sequence of
ordinals in $C^{n+1}_\delta$ such that \(\alpha = \sup_{i < \zeta}
\alpha_i\) is a limit ordinal and \(\alpha \leq \min (C^n_\delta \but
\alpha_0)\).  We shall verify that \(\alpha \in C^{n+1}_\delta \but
S^\ast_{\alpha (\ast)}\).  Let $\beta$ be the least ordinal in
\((\nacc C^n_\delta) \but (S^\ast_{\alpha (\ast)} \cut E_{n+1})\) not
less than $\alpha$.  Without loss of generality we may assume that
\(\set{\alpha_i}{0 < i < \zeta} = C^{n+1}_\delta \cut (\alpha_0,
\alpha)\).  Then
\[\set{\alpha_i}{0 < i < \zeta} = \drop (u_\beta, E_{n+1}) \cut
(\alpha_0, \alpha)\] by (\ref{gamma}) and the fact that \(\alpha \in
E_{n+1}\) and \(\beta \in C^n_\delta\).

First suppose that \(\alpha \notin C^n_\delta\).  Then \(\alpha \in
\acc u_\beta \cut \acc E_{n+1}\) which gives us \(\alpha \in
C^{n+1}_\delta\).  If \(\beta \in S^\ast_{\alpha (\ast)}\) then
\(\beta \notin E_{n+1}\) and it follows that \(\gamma^n_\delta (\beta)
\geq \alpha\) which contradicts the fact that \(C^{n+1}_\delta \cut
(\alpha_0, \alpha) \neq \emptyset.\) Thus \(\beta \notin
S^\ast_{\alpha (\ast)}\) which implies that \(u_\beta \cut
S^\ast_{\alpha (\ast)} = \emptyset\) and thereby that \(\alpha \notin
S^\ast_{\alpha (\ast)}\).  In the other case where we have \(\alpha
\in C^n_\delta\) we only need to check that \(\alpha \notin
S^\ast_{\alpha (\ast)}\).  But this is almost immediate since if
\(\alpha \in S^\ast_{\alpha (\ast)}\) we must have \(\beta > \alpha\)
which again implies the contradictory inequality \(\gamma^n_\delta
(\beta) \geq \alpha\).

For condition~\ref{CGcoh} in the case \(n = 0\) we use \ref{P2} and
note that \(\delta' \in E_0\) and \ref{P3} gives \(\drop (u_\delta,
E_0) \cut \delta' = \drop (u_\delta \cut \delta', E_0)\).  In the
induction step \(\delta' \in E_{n+1} \cut u_\delta\) implies \(\delta'
\in C^0_\delta \sub C^n_\delta\) by \ref{P2} and \ref{P3}.  Thus
\(\gamma^n_\delta (\beta) \geq \delta'\) for every \(\beta > \delta'\)
which clearly suffices.

It is also straightforward to see that \ref{CGstat} will hold when we
reach a point where no club \(E_{n+1} \sub \acc E_n\) satisfies
(\ref{nextE}).  We shall now derive a contradiction from the
assumption that (\ref{nextE}) holds for every \(n < \omega\).  Let
\(E^\omega = \Cut_{n < \omega} E_n\) and pick \[\delta \in \acc^+
(E^\omega \cut S^\ast_{\alpha (\ast)}) \cut S^\ast_{\beta (\ast)} \cut
S^\kappa_{\omega_1}.\] Let \(\gamma_n = \sup (E_{n+1} \cut \nacc
C^n_\delta \cut S^\ast_{\alpha (\ast)})\) and \(\gamma = \sup_{n <
\omega} \gamma_n\).  Because \(\delta \in S^\kappa_{\omega_1}\) we
have \(\cf \delta > \omega\) and thus by (\ref{nextE}) and the fact
that \(\delta \in E^\omega \cut S^\ast_{\beta (\ast)}\) we have
\(\gamma < \delta\).  Pick \(\alpha \in E^\omega \cut S^\ast_{\alpha
(\ast)}\) such that \(\gamma < \alpha < \delta\) and let \(\beta_n =
\min (C^n_\delta \but \alpha)\) for every \(n < \omega\).  Clearly
\(\alpha \notin \nacc C^n_\delta\) and by \ref{CGnacc} it then follows
that \(\alpha \notin C^n_\delta\).  Thus \(\beta_n > \alpha\).
Because \(\beta_n > \gamma\) we have \(\beta_n \notin E_{n+1} \cut
S^\ast_{\alpha (\ast)}\) and by (\ref{nextC}) and (\ref{gamma}) it
then follows that \(\beta_{n+1} < \beta_n\).  This is a contradiction
since \(n < \omega\) was arbitrary. \end {proof}

Fix a sequence \(\seq{C_\delta}{\delta \in S^\ast_{\beta (\ast)} \cut
\acc E^\ast}\) that satisfies the conditions of the lemma above.  Let
$R_0$ be the tree consisting of all closed bounded subsets of $\kappa$
ordered by end extension and consider the subtree
\[R = \set{c \in R_0}{\delta > \sup (c \cut \nacc C_\delta \cut
S^\ast_{\alpha (\ast)}) \mforall \delta \in S^\ast_{\beta (\ast)} \cut
\acc E^\ast}.\] Note that intersecting with \(S^\ast_{\alpha (\ast)}\)
is not essential in the definition of $R$.  As far as the argument
that follows is concerned, \(S^\ast_{\alpha (\ast)}\) could be dropped
from the definition, or more exactly, replaced by any set that
contains \(S^\ast_{\alpha (\ast)}\).  Condition~\ref{CGnacc} is
essential however.  We shall show that $R$ is a semi-proper
\((\kappa^+, \kappa)\)-tree.  We start by noting that $R$ can not have
$\kappa$-branches by condition~\ref{CGstat}.  Also, for every \(c \in
R\) and \(\alpha < \kappa\) there exists a condition \(d \in R\) such
that \(d \leq c\) and \(\max d > \alpha\).  If $R$ does not collapse
$\kappa$, it then follows that forcing with $R$ adds a
$\kappa$-branch.  We finish the proof of Proposition~\ref{part2} by
showing that $R$ does not kill stationary sets.

Let $S$ be an arbitrary stationary subset of $\kappa$, let $\dot{C}$
be an $R$-name for a club, and let \(c \in R\).  We shall find a
condition \(c^+ \leq c\) such that \(c^+ \forces \dot{C} \cut
\check{S} \neq \emptyset\).

Fix an increasing continuous sequence \(\seq{M_\eta}{\eta < \kappa}\)
of elementary submodels of $H_\chi$, where $\chi$ is some large enough
regular cardinal, such that \(\card{M_\eta} < \kappa\),
\begin {equation} \label {kinS}
M_{\eta+1} \cut \kappa \in S^\ast_{\alpha (\ast)},
\end {equation}
and \(\seq{M_\nu}{\nu \leq \eta} \in M_{\eta+1}\) for all \(\eta <
\kappa\), and $S$, $R$, $\dot{C}$, \(\alpha (\ast)\), \(\beta
(\ast)\), and the sequences \(\seq{u_\beta}{\beta < \kappa}\) and
\(\seq{C_\delta}{\delta \in S^\ast_{\beta (\ast)} \cut E^\ast}\) are
elements of $M_0$.  Pick a limit ordinal \(\delta (\ast) \in S \cut
\acc E^\ast\) such that \(M_{\delta (\ast)} \cut \kappa = \delta
(\ast)\).

The rest of the proof is divided into two cases.  In the first case we
assume that \(\delta (\ast) \notin S^\ast_{\beta (\ast)}\).  By
\ref{P2} and (\ref{Sstar}) it follows from this assumption that
\begin {equation} \label {climber}
u_{\delta (\ast)} \cut S^\ast_{\beta (\ast)} = \emptyset.
\end {equation}
We shall define a decreasing sequence \(\seq{c_i}{i < \zeta}\) of
conditions in $R$ simultaneously with an increasing sequence
\(\seq{\alpha_i}{i < \zeta}\) of ordinals such that \(c_0 = c\),
\(\sup_{i < \zeta} \alpha_i = \delta (\ast)\) and the following
conditions hold for every \(i < \zeta\):
\begin {conditionseq}
\item \label {CIM} \(c_i \in M_{\delta (\ast)}\) and \(\alpha_i <
\delta (\ast)\)
\item \label {CIF} \(\alpha_{i+1} \geq \max c_i\) and \(c_{i+1}
\forces \alpha_{i+1} \in \dot{C}\)
\item \label {CIU} \(\max c_{i+1} > \min (u_{\delta (\ast)} \but
\alpha_i)\).
\end {conditionseq}
We shall also assume that all the choices done during the construction
are made using a choice function that is in \(M_{\delta (\ast)}\).
The length $\zeta$ of the sequence will be determined during the
construction.  The successor steps in the construction are straight
forward and present no problems.

Now suppose that we are about to pick $c_i$ and $\alpha_i$ where $i$
is a limit ordinal.  Let \(\gamma = \sup_{j < i} \max c_j\).  If
\(\gamma = \delta (\ast)\) we put \(\zeta = i\) and the construction
is successfully completed.  Thus assume that \(\gamma < \delta
(\ast)\).  Clearly the only things we have to show now is that
\begin {equation} \label {inR}
\Union_{j < i} c_j \union \sing{\gamma} \in R
\end {equation}
and \(\seq{c_j}{j < i} \in M_{\delta (\ast)}\).  By
condition~\ref{CIU} \(\gamma \in u_{\delta (\ast)}\) which by
(\ref{climber}) implies that \(\gamma \notin S^\ast_{\beta (\ast)}\)
and this takes care of (\ref{inR}).  Because the sequence
\(\seq{u_\beta}{\beta < \kappa}\) is in \(M_{\delta (\ast)}\) we also
have \(u_\gamma \in M_{\delta (\ast)}\).  But since \(u_\gamma =
u_{\delta (\ast)} \cut \gamma\) and the choice function being used is
in \(M_{\delta (\ast)}\), we could obtain the same sequences
\(\seq{c_j}{j < i}\) and \(\seq{\alpha_j}{j < i}\) arguing in
\(M_{\delta (\ast)}\), if we replace \(u_{\delta (\ast)}\) by
$u_\gamma$ in condition~\ref{CIU}.  Thus \(\seq{c_j}{j < i} \in
M_{\delta (\ast)}\).  Having completed the construction we just need
to put \(c^+ = \Union_{i < \zeta} c_i \union \sing{\delta (\ast)}\)
and note that \(c^+ \forces \delta (\ast) \in \dot{C}\).

We shall now deal with the other case where we have \(\delta (\ast)
\in S^\ast_{\beta (\ast)}\).  We shall reconstruct the sequences
\(\seq{c_i}{i < \zeta}\) and \(\seq{\alpha_i}{i < \zeta}\) in a
slightly different way.  We keep conditions \ref{CIM} and \ref{CIF}
but replace \ref{CIU} by the conditions
\begin {conditionseq}
\item \label {CIUE} \(\max c_{i+1} > \max \sing{\min (u_{\delta
(\ast)} \but \alpha_i), \min (E^\ast \but \alpha_i)}\)
\item \label {CIN} \(c_i \cut \nacc C_{\delta (\ast)} = c_0 \cut \nacc
C_{\delta (\ast)}\)
\end {conditionseq}
and require that \(\alpha_0 \geq \beta (\ast)\).  We shall first deal
with the successor step since now it requires some work.  Suppose that
$c_i$ and $\alpha_i$ are defined.  Let $\eta$ be the least ordinal in
\(\delta (\ast)\) such that $c_i$ and the ordinal \(\max \sing{\min
(u_{\delta (\ast)} \but \alpha_i), \min (E^\ast \but \alpha_i)}\) are
elements of $M_\eta$ and let \(\gamma = \sup (\nacc C_{\delta (\ast)}
\cut M_{\eta + 1})\).  By \ref{CGclub}, \ref{CGnacc}, and (\ref{kinS})
\(\gamma \in \kappa \cut M_{\eta + 1}\).  Then pick $c_{i+1}$ and
$\alpha_{i+1}$ in \(M_{\eta + 1}\) such that \(c_{i+1} \leq c_i \union
\sing{\gamma + 1}\) and conditions \ref{CIF} and \ref{CIUE} are
satisfied.  In this way \(c_{i+1} \cut \nacc C_{\delta (\ast)} = c_i
\cut \nacc C_{\delta (\ast)}\) which takes care of \ref{CIN}.

Suppose then that $i$ is a limit ordinal and \(\gamma = \sup_{j < i}
\max c_j < \delta (\ast)\).  Because \(\gamma \in u_{\delta (\ast)}\)
and \(\gamma > \beta (\ast)\) we have \(\gamma \in S^\ast_{\beta
(\ast)}\) by the assumption \(\delta (\ast) \in S^\ast_{\beta
(\ast)}\).  Furthermore \(\gamma \in \acc E^\ast\) and therefore
\(C_\gamma = C_{\delta (\ast)} \cut \gamma\) by \ref{CGcoh}.  From now
on the argument is very similar to the limit step in the case \(\delta
(\ast) \notin S^\ast_{\beta (\ast)}\).  One difference is that
\(C_{\delta (\ast)}\) and $C_\gamma$ now play the role of \(u_{\delta
(\ast)}\) and $u_\gamma$ in the previous argument.  We also have to
note that the required initial segment of the sequence \(\seq{M_i}{i <
\delta (\ast)}\) is in \(M_{\delta (\ast)}\).  Of course
(\ref{climber}) does not hold now but instead condition \ref{CIN} is
designed to make (\ref{inR}) come true.  This also applies on the
final limit step where we again put \(c^+ = \Union_{i < \zeta} c_i
\union \sing{\delta (\ast)}\).  We have found the required condition
$c^+$ which concludes the proof of Proposition~\ref{part2}.

\bibliographystyle{plain}

\begin {flushleft}
Department of Mathematics \\
University of Helsinki \\
00014 Helsinki, Finland \\
{\em E-mail:} hellsten@cc.Helsinki.FI
\end {flushleft}

\begin {flushleft}
Department of Mathematics \\
University of Helsinki \\
00014 Helsinki, Finland \\
{\em E-mail:} thyttine@cc.Helsinki.FI
\end {flushleft}

\begin {flushleft}
Institute of Mathematics \\
The Hebrew University \\
Jerusalem, Israel \\
{\em and} \\
Rutgers University \\
Department of Mathematics \\
New Brunswick, NJ USA \\
{\em E-mail:} shelah@math.huji.Ac.IL
\end {flushleft}
\end {document}